
\input amstex

\documentstyle{amsppt}
\topmatter
\title ``Affine'' Hecke algebras associated to Kac-Moody groups.
    \endtitle
\author H. Garland and I. Grojnowski \endauthor
\address  Yale University, New Haven, CT 06520
\endaddress
\NoBlackBoxes
\date  March 1995 q-alg/9508019 \enddate
\NoBlackBoxes
\endtopmatter
\document

\define\bz{\bold Z}

\define\bp{\bold P}
\define\bh{\bold H}

\define\bc{\bold C}

\define\aaa{\Cal A}
\define\flag{\Cal B}

\define\aff{\widehat{H}}
\define\ring{\bz[q,q^{-1}]}

\define\barC{\bar{C}}
\define\iso{\simeq}

\define\diff{\Cal D}
\define\CO{\Cal O}
\define\CR{\Cal R}
\define\Coh{{\text{Coh}}}

\define\cstar{{\bold C^*}}
\define\gr{\text{gr}}

\define\topf{\omega}

\define\lotimes{{\otimes^L}}

\define\qbinom#1#2{\thickfracwithdelims[]\thickness0#1#2}

In this paper we give a geometric construction of Cherednik's
double affine Hecke algebra \cite{C}. This construction is suitable
for producing a certain class of modules,
and in a paper in preparation we will use the results of this paper
to describe these modules.

The construction is modeled on the fantastic construction of the usual
affine Hecke algebra as the equivariant $K$-theory of the Steinberg variety
\cite{KL1}, i.e. the Lagrangian subvariety $\Lambda$ of $T^*(\flag\times\flag)$
given by the union of the conormal bundles to the $G$-orbits on
$\flag\times\flag$.
As such, it is ``morally obvious'' what to do here---one just replaces $G$
and $\flag$ by the corresponding affine group and its flag variety.

At this level, the idea that a result similar to our theorem 4.5 must be true
was obvious to
a lot of people. In particular, Cherednik, Ginzburg, Kazhdan and Lusztig have
told us they
knew this, and we believe that Kashiwara and Tanisaki must also have known.

Thus the main content is to actually {\it do} this for the infinite dimensional
$G$ and $\flag$.
Here again, we can claim no great originality. The beautiful paper \cite{KT}
gives a lovely meaning to
$\diff$-modules on $\flag$. Using this, our main contribution is to define what
coherent sheaves
on $T^*\flag$ are. We do this using the technique of \cite{Gr}.

The fundamental importance of Cherednik's algebra suggests that it is worth
publishing this note.
In particular, the lovely structure of the modules we have found, as well as
the conjectures (work in
progress) based on this perhaps justify it. (Not to mention the remarkable work
of Cherednik!)

We would like to mention just one observation, connected with several projects
of ours. That is, that from
the point of view of this paper, Cherednik's `Fourier transform' is highly
non-trivial.

\head
1. Notation
\endhead

We will keep notation as in \cite{KT}, and suggest the reader have a copy at
hand.

Throughout this paper, $\frak g$ will be a Kac-Moody algebra with symmetrisable
generalised Cartan matrix
$(a_{ij})_{i,j\in I}$, and $G$ will be an associated group as in \cite{K}. We
will write $T$ for the
maximal torus of $G$, $X$ for its lattice of characters, and $\flag = G/B$ for
the flag variety. All
other notation will be as in \cite{KT}, and the definition of these objects are
those of \cite{K,KT}.

To this data of $(T,\frak g)$ we can associate two algebras: the Hecke algebra
$H$ and the ``affine''
Hecke algebra $\aff$. Recall that $H$ is the free $\ring$-module with basis
$T_w$, $w\in W$,
and multiplication defined by
$$ (T_s+1)(T_s-q) = 0, \qquad l(s) =1 $$
$$ T_sT_w = T_{sw}, \qquad l(sw)=l(s)+l(w).$$
Also recall that $\aff = H \otimes_{\bz}\bc[X]$, and that this is given the
unique structure of an
algebra such that the Hecke algebra $H$ and the group algebra $\bc[X]$ are
subalgebras, and
$$ fT_s- T_s{}^s\!f = (q-1) \frac{f-{}^s\!f}{1-e^{-\alpha}}$$
where $s=s_\alpha$ is a simple reflection through $\alpha$ and $e^{\alpha}$ is
the corresponding
element of $\bc[X]$.

In the case $G$ is of finite type, this is the definition of $\aff$
due to Bernstein and Zelevinski.  Cherednik's double affine Hecke
algebra \cite{C} is also of this type. It occurs when $G$ is an affine
Kac-Moody group,
with a torus of dimension $l + 1$ which contains the ``degree operator'' but
not the center.
Here $l$ is the rank of the finite reductive group $\bar{G}$, and $G$ is the
loop group of $\bar{G}$.

In particular, for $G$ of affine type with torus of dimension $l+2$,
the algebra we have just defined is a generalisation of Cherednik's
algebra to arbitrary central charge. This seems to be new.

\head
2. The Hecke algebra $H$
\endhead

This section is a minor addition to \cite{KT}.

Let $P_i$, $i\in I$, be the parabolic subgroup of $G$, defined as in \cite{KT},
and let ${}_i\Cal P = G/P_i$ be the associated generalised flag variety,
defined as in \cite{K}.

The we have a map of schemes, $r_i : \flag \to {}_i\Cal P$, which is a $\bp^1$
fibration.
The precise statement we need is

\proclaim{Lemma 2.1}
Let $F$ be a finite subset of $W$ such that $s_iF\subseteq F$, $Y=\flag^F$ is
open,
$Z = \flag_F$ is closed and $Z\subseteq Y$.
Define ${}_iY=r_i(Y)$. Then for sufficently large $k$

i) The group $U_k^{-}$ acts on ${}_iY $ locally freely, and the quotient is a
finite dimensional
smooth variety.

ii) The induced morphism $U_k^-\backslash Y \to U_k^-\backslash {}_iY $is a
$\bp^1$-bundle fibration.

iii) The natural map $r_i(Z) \to U_k^-\backslash {}_iY$ is a closed immersion.

iv) The diagram
$$\CD Y_a @>>> Y_b \\
@VVV @VVV \\
{}_iY_a @>>> {}_iY_b \endCD $$
is cartesian.
\endproclaim
We write ${}_iY_a = U_a^-\backslash {}_iY$, $Y_a = U_a^-\backslash Y$.
Also denote the natural map of (ii) as ${}_ir_l : Y_l \to {}_iY_l$, and write
(as in \cite{KT})
$p_a: {}_iY \to {}_iY_a$, $p^a_b : {}_iY_a \to {}_iY_b$, $a\geq b\geq k$ for
the natural morphisms.
These are affine.
We call such data $(Z,Y,k)$ a $P_i$-stable admissible triple.

Now, let $(Z,Y,k)$ be a $P_i$-stable admissible triple, and $\bh(Z,Y)$ be the
semisimple category of
pure mixed Hodge (right) modules, defined as in \cite{KT,2.2.11} for $\lambda =
0$.
(In other words, use pure mixed Hodge modules \cite{S} rather than
$\diff$-modules, and write
$\bh(Z,Y)$ for what is denoted $\bh(0,Z,Y)$ in \cite{KT,2.2.11}).

Let $\aaa_l$ be a pure mixed Hodge module on $Y_l$ with support contained in
$Z$.
Then define $C_i'\aaa_l$ as the sum
$$ \oplus_j {}^p\!H^j( ({}_ir_l)^*({}_ir_l)_* \aaa_l ) $$
where $({}_ir_l)^*$, $({}_ir_l)_*$ denote the usual derived functors of mixed
Hodge modules as in \cite{S}.
Note that by the decomposition theorem, $ ({}_ir_l)^*({}_ir_l)_* \aaa_l $ is
semisimple and pure of weight $\alpha$
if $\aaa_l$ is pure of weight $\alpha$, and so $C_i'\aaa_l$ is a split
semisimple pure MHM.
(Be warned that ${}^pH^0f_*$ is what is denoted $\int_f$ in \cite{KT}).

\proclaim{Lemma 2.2} If $(Z,Y)$ is $P_i$-stable admissible, then if
$\aaa\in\bh(Z,Y)$
we have $C_i'\aaa \in \bh(Z,Y)$ also.
\endproclaim

Here, if $\aaa=(\aaa_a,\gamma^a_b)$, then $C_i'\aaa =
(C_i'\aaa_a,\tilde{\gamma}^a_b)$,
and $\tilde{\gamma}^a_b$ are induced from $\gamma_b^a$ and the base change
isomorphism applied to the Cartesian diagram of lemma 2.1. The proof is
immediate (using the fact that the horizontal maps are affine space fibrations,
and so
$\int_{p^a_b} = (p^a_b)_*$).

Write $\bh=\lim_{@>Z>>}\lim_{@<Y<<} \bh(Z,Y)$ as in \cite{KT,2.2}.
Then $\bh$ is a semisimple abelian category, where every object has finite
length, and $\bh$ is
a free $\ring$ module with basis $C_w'$, $w\in W$. Here, $q$ acts by Tate
twist, and $C_w'$ represents the pure simple Hodge module which is the perverse
extension of the constant local system on $\flag_w$ with weight $l(w)$.

\proclaim{Proposition 2.3} The operators $C_i' : \bh \to \bh$ give $\bh$ the
structure of the regular $H$-module,
where $C_i'$ acts as multiplication by $T_{s_i}+1$ on $\bh$.
\endproclaim
The proof is standard and omitted. (It consists of enlarging $\bh$ to allow
non-semisimple perverse sheaves,
so that the constant local systems on $\flag_w$ extended by $0$, call them
$T_w$, are in this new category $\tilde{\bh}$.
Then the Grothendieck group of $\tilde{\bh}$ is the same as that of $\bh$, and
the elements $T_w$ clearly satisfy
$T_sT_w = T_{sw}$ if lengths add, and $(T_s+1)^2$ is as claimed by standard
properties of $\bp^1$ fibrations.)

\head
3. Definition of $\Coh^{T\times \cstar}(\Lambda)$
\endhead

3.1

Let $X$ be a variety on which our torus $T$ acts. A $T$-equivariant coherent
sheaf on $X$ is a
coherent sheaf on the quotient stack $X/T$; explicitly it is a pair
$(\aaa,\phi)$ where $\aaa$ is a coherent sheaf
on $X$ and $\phi:p^*\aaa\to a^*\aaa$ is an isomorphism satisfying the usual
compatibilities.
Here $p:T\times X \to X$, $(t,x)\mapsto x$ and
 $a:T\times X \to X$, $(t,x)\mapsto tx$ are the projection and action maps.

Denote the category of $T$-equivariant coherent sheaves $\Coh^T(X)$. This is an
exact category,
write $D^T(\CO_X)$ for its derived category and $K^T(X)$ for the Grothendieck
group of $\Coh^T(X)$
which is canonically isomorphic to the Grothendieck group of $D^T(\CO_X)$.

If $f:X\to Z$ is a morphism of $T$-spaces, we denote by $f^*$ the left derived
functor of the
usual pullback of $\CO$-modules, and by $f_*$ the right derived functor of the
usual pushforward of $\CO$-modules.
Recall that if $f$ is a regular morphism, $f^*\aaa$ is a finite complex of
coherent $\CO$-modules if $\aaa$ is,
and so $f^*: D^T(\CO_Z) \to D^T(\CO_X)$ in this case, and if the restriction of
$f$ to the support of $\aaa$ is
proper, then $f_*\aaa$ is a finite complex of coherent $\CO$-modules if $\aaa$
is, so
$f_* : D^T(\CO_X)\to D^T(\CO_Z)$ in this case. We refer to \cite{T} for
details.

Now let $Y\subseteq X$ be a closed $T$-stable subvariety. Write $D^T(\CO_X, Y)$
for the triangulated
subcategory of $D^T(\CO_X)$ consisting of sheaves whose cohomology has support
in $Y$, and write
$K^T(X,Y)$ for its Grothendieck group.

If $i:Y\hookrightarrow X$ denotes the inclusion, then $i_*: D^T(\CO_Y) \to
D^T(\CO_X,Y)$
induces an isomorphism $i_*: K^T(Y)\to K^T(X,Y)$. In particular, $K^T(X,Y)$
does not depend on the
ambient variety $X$.
(We recall that a sheaf on $X$ with support on $Y$ need not be of the form
$i_*\aaa$, and so there is a small
subtlety in these statements).

3.2

Let $(Z,Y,k)$ be an admissible triple. Then $Z=\flag_F$, for $F$ a finite
subset of $W$. Define, for $l \geq k$
and $F'\subseteq F$
$$ \Lambda^l_{F'} = \bigcup_{w\in F'} T^*_{\flag_w}Y_l \subseteq T^*Y_l. $$
Write $\Lambda_Z^l = \Lambda_F^l$, and note that the group $T\times\cstar$ acts
on $\Lambda_{F'}^l$.

Now, the affine fibration $p^a_b : Y_a \to Y_b $ gives rise to a correspondence
$$ T^*Y_a @<{\alpha^a_b}<< Y_a \times_{Y_b} T^*Y_b @>{\beta^a_b}>> T^*Y_b $$
where we suppose $a \geq b \geq k$.

Ignore the $T$-action for a moment. Then, choosing a small open set
$U_b\subseteq Y_b$,
we get an open set $U_a = (p^a_b)^{-1}(U_b)$ of $Y_a$, and we can identify
$$ U_a \iso U_b \times \bc^n$$
where $p_b^a : U_a \to U_b$ is projection onto the second factor, and $Z\cap
U_a \hookrightarrow U_a$ is
the map $z \mapsto (z,0)$.

It follows that $U_a \times_{U_b} T^* U_b \iso T^*U_b \times \bc^n \times 0$,
$T^*U_a \iso  T^*U_b \times \bc^n \times \bc^n$, and that $\alpha^a_b$ is the
obvious embedding,
$\beta^a_b$ the evident projection, and that
$\Lambda^a_Z \cap T^*U_a \iso (\Lambda_Z^b \cap T^*U_b) \times 0 \times \bc^n$.

In particular, we get
$(\alpha^a_b)^{-1} (\Lambda^a_Z) \iso \Lambda^b_Z \times 0 \times 0$, and
$\beta^a_b$ induces
an isomorphism of this onto $\Lambda^b_Z$.

Thus, if we denote by $I^a_b : D(\CO_{T^*Y_a}, \Lambda^a_Z) \to D(\CO_{T^*Y_b},
\Lambda^b_Z)$
the map $(\beta^a_b)_*(\alpha^a_b)^*$, then this does take a coherent sheaf on
$T^*Y_a$ supported on
$\Lambda^a_Z$ to a complex of sheaves on $T^*Y_b$, with cohomology supported on
$\Lambda_Z^b$.

Further, the induced map $I_b^a : K(T^*Y_a,\Lambda^a_Z) \to
K(T^*Y_b,\Lambda^b_Z)$ is an
isomorphism (as restricting to the zero section of a vector bundle induces an
isomorphism in $K$-theory).

Now, recall that we have a $T$-action. We still get maps
$$ K^{T\times \cstar}(\Lambda^a_Z) @>{\sim}>>
K^{T\times\cstar}(T^*Y_a,\Lambda^a_Z) @>{I^a_b}>>
K^{T\times\cstar}(T^*Y_b,\Lambda^b_Z) @<{\sim}<< K^{T\times
\cstar}(\Lambda^b_Z)$$
and in fact

\proclaim{Lemma 3.3} $I^a_b :  K^{T\times\cstar}(T^*Y_a,\Lambda^a_Z) \to
K^{T\times\cstar}(T^*Y_b,\Lambda^b_Z) $
is an isomorphism.
\endproclaim
But this follows from the previous discussion, either by restricting to
$T$-stable neighbourhoods
$U_b$  of $Y_b$ and observing we still have a splitting (Sumihoro's theorem),
or by restricting to $T$-fixpoints and using the embedding
$K^T(X)\hookrightarrow K^T(X^T)$ and the preceding arguments.

Now we can define a group depending on $(Y,Z,k)$ as the limit of these
isomorphisms; i.e. an element here is a
sequence $x_a \in K^{T\times\cstar}(T^*Y_a,\Lambda^a_Z)$, $a \geq k$ such that
$I^a_bx_a=x_b$, $a\geq b$.
We can take the limit over $k$. Further, if $Y'$ is an admissible open set
containing $Y$, then the group
defined for $(Y,Z)$ is canonically isomorphic to that for $(Y',Z)$. So we can
take the limit over admissible open sets $Y$
containing $Z$. Call this resulting group
$$ K^{T\times\cstar}(\Lambda_Z). $$
Finally, if $Z'$ is an admissible closed set containing $Z$, the the closed
embeddings
 $i_a:\Lambda^a_Z\hookrightarrow\Lambda^a_{Z'}$ give rise to an embedding
 $i_* :  K^{T\times\cstar}(\Lambda_Z) \hookrightarrow
K^{T\times\cstar}(\Lambda_{Z'}) $.
Take the limit over admissible closed sets $Z$ to get (compare \cite{KT,2.2})
$$   K^{T\times\cstar}(\Lambda) = \lim_{@>>Z>} K^{T\times\cstar}(\Lambda_Z). $$

3.4
This Grothendieck group is somewhat crude. However, the same procedure works
nicely to define $\Coh^{T\times\cstar}(\Lambda)$,
``actual'' coherent sheaves with finite support on $\flag$ and finite cosupport
along the cotangent directions in
 $T^*\flag$ (and supported in $\Lambda$).

The point is to define for each admissible triple $(Y,Z,k)$ a ``coherent
sheaf'' $\aaa$ as a family
$\aaa= ((\aaa_l)_{l\geq k}, (\gamma^a_b)_{a\geq b\geq k})$, where $\aaa_l$ is a
coherent sheaf on $T^*Y_l$ supported
on $\Lambda_Z^l$ and
$$ \gamma^a_b : I^a_b\aaa_a @>\sim>> \aaa_b $$
is an isomorphism satisfying the chain condition $\gamma^a_cI^a_b = \gamma^b_c
I^b_c \gamma^a_b$.
Note that in order for $\gamma^a_b$ to exist, $\aaa_a$ must be smooth (locally
constant) in the directions transverse
to the bundle embedding $Y_a \times_{Y_b} T^*Y_b \hookrightarrow T^*Y_a$,
otherwise $I^a_b\aaa_a$ will be a complex of
sheaves.

Nonetheless, enough such families exist; one can define morphisms in the
obvious way, and take limits as above, getting
the abelian category $\Coh^{T\times\cstar}(\Lambda)$. Then the Grothendieck
group of $\Coh^{T\times\cstar}(\Lambda)$
is $K^{T\times\cstar}(\Lambda)$, as defined above.

We will not need this precision, but it's nice to know its there.

\head
4. Action of $\aff$ on $K^{T\times\cstar}(\Lambda)$.
\endhead

4.1

We now define an action of the Hecke algebra $H$ on
$K^{T\times\cstar}(\Lambda)$.

Let $(Z,Y,k)$ be a $P_i$-stable admissible triple, $i\in I$. Let $a\geq b\geq
k$.
Define
$$ {}_i\CR_a = \{ (x,y) \in Y_a \times Y_a \mid {}_ir_a(x)={}_ir_a(y) \},$$
and let $\pi_a, \pi'_a$ denote the first and second projections. Then by lemma
2.1, $\pi_a : {}_i\CR_a \to Y_a$
is a $\bp^1$-bundle. We obtain an induced correspondence on cotangent bundles
$$ T^*Y_a @<{\pi_a}<< T^*_{{}_i\CR_a}(Y_a\times Y_a) @>{\pi'_a}>> T^*Y_a.$$
One can easily check that $\pi'_a(\pi_a^{-1}(\Lambda^a_Z)) \subseteq
\Lambda^a_Z$,
and so we get a map
$\barC_i' : D^{T\times\cstar}(\CO_{T^*Y_a},\Lambda^a_Z) \to
D^{T\times\cstar}(\CO_{T^*Y_a},\Lambda^a_Z)$
defined by
$$ \barC'_i\aaa = (\pi_a')_* (\pi^*_a\aaa \lotimes \rho^{-1}
\topf_{{}_i\CR_a/Y_a}) [1] $$
where $\topf_{{}_i\CR_a/Y_a} $ is the bundle of relative top forms
$\topf_{{}_i\CR_a/Y_a} = \pi_a^*(\topf_{Y_a}^{\otimes^{-1}})
\otimes_{\pi_a^{-1}\CO_{Y_a}} \topf_{{}_r\CR_a}$,
pulled back via $\rho : T^*_{{}_i\CR_a}(Y_a\times Y_a) \to {}_i\CR_a$.
(Note that here we use the first projection, unlike in \cite{Gr,2}, as we are
using right $\diff$-modules).

Furthermore, one can easily check that if $\aaa\in
D^{T\times\cstar}(\CO_{T^*Y_a},\Lambda^a_Z)$, $a \geq b \geq k$,
the $\barC'_i I^a_b\aaa = I^a_b\barC_i'\aaa$; the diagram chase is omitted.

It follows we have well defined maps $\barC_i' : K^{T\times\cstar}(\Lambda) \to
 K^{T\times\cstar}(\Lambda)$

4.2

Now let $w\in W$ be such that $\flag_w \subseteq Z$. The pure mixed Hodge
module which is the perverse extension
of the constant local system on $\flag_w$, with weight $l(w)$, is part of a
compatible family of such,
which we denoted $C'_w$ in \S2.

If $(M,F,\dots)$ represents this in $Y_a$, then $\gr_FM\in
\Coh^{T\times\cstar}(\Lambda_Z^a)$, and by definition
of $\Coh^{T\times\cstar}(\Lambda)$ and \cite{Gr,2.1.2} these individual
coherent sheaves patch togethor to to give an element
of $\Coh^{T\times\cstar}(\Lambda)$. Denote this $\gr\, C_w'=\barC_w'$.

It is clear that $\{\barC'_w \mid w\in W\}$ are $\ring$-linearly independant in
$K^{T\times\cstar}(\Lambda)$,
and so we have defined an embedding $\gr : \bh \to K^{T\times\cstar}(\Lambda)$.

Moreover,  by the compatibility of $\gr$ and correspondences (see \cite{Gr,2}),
the action of the Hecke algebra $H$
by right multiplication (by $C'_i$) induces an action on $\gr\bh$; this is
precisely the action by the
operators $\barC_i'$ defined above.

4.3

Let $\lambda:T\to \cstar$ be a character of $T$, $\CO_{Y_l}(\lambda)$ the
invertible line bundle on $Y_l$
defined in \cite{KT,2.2.5}. Then if $p: T^*Y_l \to Y_l$ denotes the canonical
projection,
we get an operator
$\Theta_\lambda : D^{T\times\cstar}(\CO_{T^*Y_a},\Lambda^a_Z) \to
D^{T\times\cstar}(\CO_{T^*Y_a},\Lambda^a_Z)$,
$ \aaa \mapsto \aaa \otimes p^*\CO_{Y_l}(\lambda)$.

Clearly $\Theta_\lambda I^a_b = I^a_b \Theta_\lambda$, and
$\Theta_\lambda\Theta_\mu = \Theta_{\lambda+\mu}$.
This gives an action of $\bc[X]$, the group algebra of $X$, on
$K^{T\times\cstar}(\Lambda)$.
Recall that $W$ also acts on $\bc[X]$.

\proclaim{Proposition 4.4} If $f \in \bc[X]$, and $s= s_\alpha$ is a simple
reflection in $W$, then as operators
on $K^{T\times\cstar}(\Lambda)$
$$ fT_s- T_s{}^s\!f = (q-1) \frac{f-{}^s\!f}{1-e^{-\alpha}},$$
where $T_s = \barC_s -1$.
\endproclaim
Let $Y @>{\pi}>> \bar{Y}$ be a $\bp^1$-fibration,
 $\CR=\{(x,y)\in Y\times Y \mid \pi(x)=\pi(y) \}$,
$T^*Y @<<< T_{\CR}^*(Y\times Y) @>>> T^*Y$ the associated correspondence. Then
the proposition is true in the
generality of such correspondences. It is the semisimple rank 1 ``affine''
version of \cite{KL1,1.3o2}
and easily follows from standard facts about the cohomology of $\bp^1$-bundles
and the Koszul complex.
We omit the short proof.

It follows that we have an action of $\aff$ on $\aff.\barC'_1 \subseteq
K^{T\times\cstar}(\Lambda)$.
\proclaim{Theorem 4.5} The map $\aff \to  K^{T\times\cstar}(\Lambda)$, $h
\mapsto h.\barC_1'$ is
an isomorphism.
\endproclaim
\demo{Proof}
Refine the Bruhat order $\leq$ on $W$ to a total order on $W$, still denoted
$\leq$. Write $x < y$ if
$x \leq y$ and $x \neq y$. Now, if $(Z,Y,k)$ is an admissible triple, $a\geq
k$, $\flag_w \subseteq Z$
then
$$ K^{T\times\cstar}(\Lambda^a_w) =  K^{T\times\cstar}(T^*_{\flag_w}Y_a) =
K^{T\times\cstar}(\flag_w) =  K^{T\times\cstar},$$
the last two equalities as $T^*_{\flag_w}Y_a$ is an affine space bundle over
$\flag_w$, and $\flag_w$ is itself affine.
It follows \cite{T} that we have a short exact sequence in $K$-homology
$$ 0 @>>>  K^{T\times\cstar}(\Lambda^a_{<w}) @>>>
K^{T\times\cstar}(\Lambda^a_{\leq w}) @>{j^*}>>  K^{T\times\cstar}(\Lambda^a_w)
@>>> 0.$$
Now, if $(M,F,\dots)$ is the mixed Hodge module corresponding to the basis
element $C'_z\in\bh$, then it is clear
that $\gr\, C'_z \in K^{T\times\cstar}(\Lambda^a_{<w}) $ if $z<w$, and
$ j^*(\Theta_\lambda\,\gr\,C'_w),  \lambda \in X$
form a basis of $K^{T\times\cstar}(\Lambda^a_w)$ over $\ring$.

As all these short exact sequences are compatible with the maps $I^a_b$, we see
that we have produced a filtration
of $K^{T\times\cstar}(\Lambda)$, and an isomorphism from $\aff \to \gr\,
K^{T\times\cstar}(\Lambda)$
(as $\Theta_\lambda C'_w$, $\lambda\in X, w\in W$ form a basis of $\aff$). The
theorem follows.
\enddemo

\head
5. Standard modules.
\endhead

Define $T^*\flag = G\times^B\frak n = \{ (x,gB)\in\frak g \times \flag \mid
g^{-1}xg \in \frak n \}$, and
$\Lambda = \{ (x,gB) \in T^*\flag \mid x \in \frak n\}$.
Then one can give $\Lambda$ the structure of a scheme as in \cite{K}.

Given $(Z,Y,k)$ an admissible triple, let $\Lambda_Y = \Lambda \cap (\frak n
\times Y)$, $l \geq k$.
Then $p_l:Y\to Y_l$ gives rise to the correspondence
$$ T^*Y @<\alpha<< Y \times_{Y_l} T^*Y_l @>\beta>> T^*Y_l$$
and one may check
\proclaim{Lemma 5.1} $\beta$ defines an isomorphism between
$\alpha^{-1}(\Lambda_Y \cap (\frak n \times\flag_w))$
and $\Lambda_w^l$, for any $\flag_w \subseteq Z$.
\endproclaim

Similarly, if $(s,q) \in T\times \cstar$, then one may define its fixpoints on
$Y_l, Y, \dots$
and if $n\in\frak n$, $s.n = qn$, then we can define $K_*(\flag^s_n)$ and an
action of $\aff$ on it,
even though the variety $\flag^s_n$ need not be rational \cite{KL2}.


In the sequel to this paper, we will study these ``standard'' modules and their
irreducible quotients,
in the case when $G$ is affine.

\enddocument

\ref\key BBD \by A. Beilinson, J. Bernstein and P. Deligne\paper
Faisceaux pervers \jour Ast\'erisque\vol 100\yr 1982\endref
\ref\key K \by M. Kashiwara\paper On crystal bases of the $q$-analogue
of universal enveloping algebras \jour Duke Math. J. \vol 63
\yr 1991\pages 465-516\endref
\ref\key L1 \by G. Lusztig \paper Character sheaves, II \jour Adv. in Math.
\yr 1985  \vol 57 \pages 226-265\endref
\ref\key L2 \by G. Lusztig \paper Cuspidal local systems and graded Hecke
algebras\jour Publications Math\'ematiques
 \vol 67\yr 1988\pages 145-202\endref
\ref\key L3 \by G. Lusztig \paper Canonical bases arising from quantized
enveloping algebras \jour J. Amer. Math. Soc.\vol 3\yr 1990\pages 447-498
\endref
\ref\key L4 \by G. Lusztig \paper Canonical bases arising from quantized
enveloping algebras, II \jour Progr. Theor. Phys. Suppl.\vol 102
\yr 1990\pages 175-201\endref
\ref\key G1\by I. Grojnowski\paper The coproduct for quantum $GL_n$
\paperinfo preprint 1992\endref
\ref\key G2\by I. Grojnowski\paper Braid group action on quantum $GL_n$
\paperinfo preprint 1991\endref

\ref\key B\by J. Beck \paper Braid group action and quantum affine algebras
\jour Communications in Math. Physics \paperinfo to appear \endref
\ref\key BB\by A. Beilinson and J. Bernstein \paper A proof of Janzten
conjectures
\jour Adv. in Soviet Math. \vol 16\yr 1993 \pages 1--50 \endref
\ref\key BL\by J. Bernstein and V. Lunts \paper Equivariant sheaves and
functors
\inbook LNM 1578 \publ Springer-Verlag \yr 1994
\endref
\ref\key DF\by J. Ding and I. Frenkel\paper Isomorphism of two realizations of
quantum affine algebra $U_q\widehat{\frak g\frak l(n)}$
\jour Communications in Math. Physics \yr 1993\vol 156\pages 277--300\endref
\ref\key D\by V. Drinfeld\paper A new realisation of Yangians and quantized
affine algebras
\jour Soviet Math. Dokl. \vol 36\yr 1988 \pages 212--216\endref
\ref\key F\by M. Finkelberg\paper Microlocalisation of linear sheaves
\paperinfo letter to V. Schectman
\yr 1993 \endref
\ref\key GV\by V. Ginzburg and E. Vasserot
\paper Langlands reciprocity for affine quantum groups of type $A_n$
\jour International Math. Research Notes \vol 3\yr 1993\pages 67--85 \endref
\ref\key G1\by I. Grojnowski\paper Representations of quantum affine algebras
\paperinfo Yale University course notes, book in preparation \yr 1994\endref
\ref\key G2\by I. Grojnowski\paper Character sheaves on symmetric spaces
\paperinfo MIT thesis \yr 1992 \endref
\ref\key GL\by I. Grojnowski and G. Lusztig\paper A Comparison of bases of
quantized enveloping algebras \jour Contemp. Math 153 \pages 11--19\yr
1993\endref
\ref\key La1\by G. Laumon \paper Sur le categorie derivee des $\diff$-modules
filtres
 \inbook LNM 1016 \publ Springer-Verlag \yr 1983 \pages 151--237\endref
\ref\key La2\by G. Laumon \paper Transformations canoniques et specialisation
pour les
$\diff$-modules filtres \jour Asterisque\vol 130\yr 1985\pages 56--129
\endref
\ref\key L1 \by G. Lusztig \paper Quivers, perverse sheaves and quantized
enveloping algebras\jour J. Amer. Math. Soc. \vol 4 \yr 1991
\pages 365-421\endref
\ref\key L2 \by G. Lusztig \book Introduction to Quantum groups
\publ Birkhauser\publaddr Boston \bookinfo Progress in Math \vol 110 \yr 1993
\endref
\ref\key L3 \by G. Lusztig \paper Character Sheaves\jour Adv. in Math. \vol 56
 \yr 1985\pages 195--237\endref
\ref\key Na\by H. Nakajima\paper Instantons on ALE spaces, quivers, and
Kac-Moody algebras
\paperinfo preprint 1994 \endref
\ref\key S1 \by M. Saito\paper Modules de Hodge polarisables \vol 24\yr 1988
\pages 849--995
\jour Publ. RIMS\endref
\ref\key V\by J. L. Verdier\paper Specialisation de faisceaux et monodromie
moderee
\pages 332--364 \vol 101\jour Asterisque \endref
\ref\key Vo\by D. Vogan\paper Irreducible characters of semisimple Lie groups
 IV\jour Duke Math Jour. \vol 49\yr 1982\pages 943--1073 \endref

\enddocument